%

\input amstex
\documentstyle{amsppt}
\magnification=\magstep1 
 
\def\binrel@#1{\setbox\z@\hbox{\thinmuskip0mu
\medmuskip\m@ne mu\thickmuskip\@ne mu$#1\m@th$}%
\setbox\@ne\hbox{\thinmuskip0mu\medmuskip\m@ne mu\thickmuskip
\@ne mu${}#1{}\m@th$}%
\setbox\tw@\hbox{\hskip\wd\@ne\hskip-\wd\z@}}
\def\overset#1\to#2{\binrel@{#2}\ifdim\wd\tw@<\z@
\mathbin{\mathop{\kern\z@#2}\limits^{#1}}\else\ifdim\wd\tw@>\z@
\mathrel{\mathop{\kern\z@#2}\limits^{#1}}\else
{\mathop{\kern\z@#2}\limits^{#1}}{}\fi\fi}
\def\underset#1\to#2{\binrel@{#2}\ifdim\wd\tw@<\z@
\mathbin{\mathop{\kern\z@#2}\limits_{#1}}\else\ifdim\wd\tw@>\z@
\mathrel{\mathop{\kern\z@#2}\limits_{#1}}\else
{\mathop{\kern\z@#2}\limits_{#1}}{}\fi\fi}
\def\circle#1{\leavevmode\setbox0=\hbox{h}\dimen@=\ht0
\advance\dimen@ by-1ex\rlap{\raise1.5\dimen@\hbox{\char'27}}#1}
\def\sqr#1#2{{\vcenter{\hrule height.#2pt
     \hbox{\vrule width.#2pt height#1pt \kern#1pt
       \vrule width.#2pt}
     \hrule height.#2pt}}}

\def\force{\hbox{$\|\hskip-2pt\hbox{--}$\hskip2pt}}  
\def\notforce{\,\nobreak\not\nobreak\!\nobreak\force}
 
\baselineskip 24pt

\define\k{\kappa}
  \define\a{\alpha}

\topmatter
\title
Full Reflection at a Measurable Cardinal\\
\endtitle
\author  Thomas Jech and Ji\v r\' i Witzany   \endauthor
\thanks The first author was supported by NSF grant number
DMS-8918299. The second author was partially supported by the first author's
NSF grant. Both authors wish to thank S. Baldwin, W. Mitchell and H. Woodin
for their comments on the subject of this paper.\endthanks
\affil     The Pennsylvania State University  and
Charles University (Prague)                   \endaffil
\address Department of Mathematics, The Pennsylvania State University,
University Park, PA 16802                       \endaddress

\abstract
A stationary subset $S$ of a regular uncountable cardinal $\k$ {\it reflects 
fully}
at regular cardinals if for every stationary set $T \subseteq \k$ of higher
order consisting of regular cardinals there exists an $\a \in T$ such that
$S \cap \a$ is a stationary subset of $\a$. {\it Full Reflection} states that
every stationary set reflects fully at regular cardinals. We will prove that 
under
a slightly weaker assumption than $\k$ having Mitchell order $\k^{++}$ it is
consistent that Full Reflection holds at every $\lambda \leq \k$
and $\k$ is measurable. 
\endabstract 
\email  jech\@math.psu.edu, witzany\@math.psu.edu                    \endemail 
\endtopmatter
 
\document

\subhead
1. Definitions and results
\endsubhead
 
It has been proved in [M82] that reflection of stationary sets is
a large cardinal property. Reflection of stationary
subsets of $\omega_n \; (n\geq2)$ and $\omega_{\omega +1}$
has been
investigated
in [M82] and [JS90] and consistency strength of Full Reflection at regular 
cardinals
at a Mahlo cardinal has been characterized in [JS92].
In this paper we address the
question of Full Reflection at a measurable cardinal.

If $S$ is a stationary subset of a regular uncountable cardinal $\k$ then
{\it the trace of} $S$ is the set

$$ Tr (S) = \{ \a < \k ; \; \; S \cap \a \; \text{is stationary
in}  \; \a \} $$
(and we say that $S$ {\it reflects at} $\a$).  If $S$ and $T$ are
both stationary,
we define
$$ S < T \;\; \text{if for almost all} \; \a \in T, \; \; \a \in
Tr(S) $$
and say that $S$ {\it reflects fully} in $T$.
(Throughout the paper, ``for almost all" means ``except for a
nonstationary set of
points").

\proclaim{Lemma 1.1} {\rm ([J84])} The relation $<$ is well founded. 
\endproclaim

\demo{Proof}
By contradiction suppose there is a sequence of stationary sets such that
$$A_1 > A_2 > A_3 > \cdots $$
It means that there are clubs $C_n$ such that
$$ A_n \cap C_n \subseteq Tr(A_{n+1}) \; \; \text{for} \; n=1,2,\dots $$
If $C \subseteq \k$ is a club let us denote $C'=\{ \a < \k ;\; C\cap \a 
\text{ is unbounded in }  \a \}$,
( $C'$ is again a club) and put
$$\tilde A_n = A_n \cap C_n \cap C_{n+1}' \cap C_{n+2}'' \cap \cdots \; \; \; 
\; \text{for} \; n=1,2, \dots$$
Then all $\tilde A_n$ are stationary. Observe that $\a \in Tr(S)$ implies 
$cf(\a)> \omega$ and 
$Tr(S\cap C)=Tr(S) \cap C'$ where $C$ is any club. Now it is easy to verify
that
$$ \tilde A_n \subseteq Tr(\tilde A_{n+1}) \;\;\;\; \text{for } n=1,2,\dots $$
Let $\a_n = \min( \tilde A_n)$. Since $\tilde A_{n+1} \cap \a_n$ is stationary 
in $\a_n$,
the ordinal $\a_{n+1}$ must be less then $\a_n$. We obtain a sequence of 
ordinals
$$ \a_1 > \a_2 > \a_3 > \cdots $$
- a contradiction. \qed
\enddemo

 {\it The order} $o(S)$ of a stationary set of regular cardinals is defined
as the rank of S in relation $<$:
$$ o(S)=\sup \{ o(T)+1 ; \; T \subseteq Reg(\k ) \text{ stationary and }  
T<S \} .$$
For a stationary set $T$ such that $T\cap Sing(\k )$ is stationary we
define $o(T)=-1$. {\it The order of } $\k$ is then defined as
$$ o(\k) = \sup\{ o(S)+1 \; ;  S \subseteq \k \; \text{is stationary} \}.$$
Note that if $Tr(S)$, where $S \subseteq Reg(\k)$, is stationary then
$o(S) < o(Tr(S))$ because $S<Tr(S)$. It follows from [J84] that the order
$o(\k)$ provides a natural generalization of the Mahlo hierarchy: $\k$ is
exactly $o(\k)$-Mahlo if $o(\k) < \k^{+}$ and  greatly Mahlo
if $o(\k) \geq \k^{+}$.

We say that a stationary set $S$ {\it reflects fully} at regular cardinals 
if for
any stationary set $T$ of regular cardinals $o(S)<o(T)$ implies $S<T$.

\proclaim{Axiom of Full Reflection at $\k$}  Every stationary subset of $\k$ 
reflects
fully at regular cardinals. \endproclaim

Following [J84] we say that  a stationary set $E$ is {\it canonical of order} 
$\nu$ if $E$ is hereditarily 
of order $\nu$ (i.e. $o(X)=\nu$ for every stationary $X \subseteq E$) and $E$ 
meets every 
stationary set of order $\nu$.

The existence of canonical stationary sets of order less then $\k ^{+}$ (if a 
set of
such order exists) is proved in [BTW76] and [J84].
In the model constructed in Section 3 we get a sequence of stationary sets
with the following properties:

\proclaim{Lemma 1.2}  Let  $\langle E_{\delta} ; \;
-1\leq \delta < \theta \rangle$ be a maximal antichain
of stationary subsets of $\lambda$ such that

(i) $E_{-1}=Sing(\lambda), \; E_\delta \subseteq Reg(\lambda)$ for $\delta
\geq 0,$

(ii) for any $\delta \geq 0$ the set $Tr(E_\delta)\cap E_\delta$ is 
nonstationary,

(iii) if $S\subseteq E_\delta$ is stationary and $-1 \leq \delta < \delta'$
then $S<E_{\delta'}.$
\newline
Then each $E_\delta$ is a canonical stationary set of order $\delta,$
$o(\lambda)=\theta$ and Full Reflection holds at $\lambda.$
\endproclaim

\demo{Proof}
We will prove the lemma in several steps. Obviously $E_\delta < E_{\delta'}$
if $\delta<\delta'.$

\proclaim{Claim 1}
Let $T\subseteq Reg(\lambda)$ be a stationary set such that
$T\cap E_{\delta'}$ is nonstationary for $\delta' \leq \delta$ and 
$S\subseteq E_\delta$ stationary. Then $S<T.$
\endproclaim
\demo{Proof}
We need to prove that $T\setminus Tr(S)$ is nonstationary.
But $(T\setminus Tr(S))\cap E_{\delta'}$ is nonstationary for $\delta' \leq 
\delta$ because $T\cap E_{\delta'}$ is nonstationary, and for $\delta'>\delta$
the set $(T\setminus Tr(S)) \cap E_{\delta'}$ is nonstationary because
$E_{\delta'} \setminus Tr(S)$ is nonstationary.
Consequently $T\setminus Tr(S)$ is nonstationary.
\enddemo

\proclaim{Claim 2}
If $S\subseteq E_\delta$ is stationary then $o(S)=o(E_\delta).$
\endproclaim
\demo{Proof}
Suppose the claim holds for $\delta'<\delta$ and that
for some $S\subseteq E_\delta$ stationary
$o(S)>o(E_\delta).$
Then there is $T<S$ such that $o(T)=o(E_\delta).$
By the induction hypothesis
$T\cap E_{\delta'}$ must be nonstationary for $\delta'< \delta$
($T_1 \subseteq T_2$ stationary implies $o(T_1)\geq o(T_2)$).
Moreover $E_\delta \cap Tr(T\cap E_\delta)$ is nonstationary. Thus
$ S\setminus Tr(T\setminus E_\delta)$ must be nonstationary because
$S\setminus Tr(T) = (S\setminus Tr(T\setminus E_\delta))\cap
(\lambda \setminus (E_\delta \cap Tr(T\cap E_\delta))$
is nonstationary. It means that $(T\setminus E_\delta)<S$
but by claim 1 $S< (T\setminus E_\delta)$ which is a contradiction
with well-foundedness of $<$.
\enddemo

\proclaim{Claim 3}
$o(E_\delta)=\delta$ for $\delta < \theta_\lambda.$
\endproclaim
\demo{proof}
Suppose by induction that $o(E_{\delta'})=\delta'$ for $\delta' < \delta .$
Certainly $o(E_\delta)\geq \delta$, suppose by contradiction that there is
a set $T<E_\delta$ such that $o(T)=\delta.$ As in the proof of claim 2
we can suppose that $T\cap E_{\delta'}$ is nonstationary for
$\delta' \leq \delta.$ But it implies by claim 1 that $E_\delta < T$
- a contradiction.
\enddemo

It follows from these claims that each $E_\delta$ is a canonical stationary
set of order $\delta.$ Any $S\subseteq Reg(\lambda)$ stationary must have a 
nonstationary intersection with some $E_\delta$ which means $o(S)\leq\delta$
and so $o(\lambda)=\theta_\lambda,$ actually 
$$o(S)=\min\{ \delta<\theta_\lambda ;\; E_\delta \cap S \text{ is stationary}\}.
$$
Finally let $S\subseteq \lambda,$ $T\subseteq Reg(\lambda)$ be stationary
and $\delta=o(S) < o(T)$ then by claim 1 $S\cap E_\delta < T$ which implies
$S<T.$
\qed
\enddemo

To state our result we need to review the definition of Mitchell order and
of a coherent sequence.

If $\;\Cal U, \; \Cal V$ are two measures on $\k$ then $\Cal U \vartriangleleft 
\Cal V$ is defined
as $\Cal U \in V^{\k}{/\Cal V}$. The transitive relation $\vartriangleleft$ 
is known to be well-founded (see [Mi74]). 
{\it The Mitchell order of} $\k$ is then defined
as the rank of this relation on measures over $\k$.

A coherent sequence of measures is a function $\overrightarrow{\Cal U}$ with
domain of the form $\{ (\alpha, \beta) ; \; \alpha < l(\Cal U) \text{ and }  
\beta < o^{\Cal U}(\alpha) \}$
for some ordinal $l(\Cal U)$ and a function $o^{\Cal U}(\cdot)$ such that
$$
\align
(i)&\text{ For all } (\alpha , \beta ) \in \text{dom} 
\overrightarrow{\Cal U} \; \; \; \Cal U^{\alpha}_{\beta} =
\overrightarrow{\Cal U}(\alpha, \beta) \text{ is a measure on }\alpha,\\
(ii)&\text{ if $j$ is the canonical embedding $j:V \rightarrow V^{\alpha}
{/\Cal U^{\alpha}_{\beta}}$
then}\\ 
&\text{ $j(\overrightarrow{\Cal U})\restriction(\alpha +1) = \overrightarrow
{\Cal U} \restriction 
(\alpha, \beta)$ where}\\
&\text{ $\overrightarrow{\Cal U}\restriction (\alpha, \beta) = \overrightarrow 
{\Cal U} \restriction
\{ (\alpha ', \beta '); \; \alpha ' < \alpha \text{ or } \alpha ' = \alpha 
\text{ and }
\beta ' < \beta \}$ and}\\
&\text{ $\overrightarrow{\Cal U} \restriction (\alpha +1) = \overrightarrow 
{\Cal U}
\restriction (\alpha +1, 0)$.}
\endalign
$$

Observe that in particular $\Cal U ^{\alpha}_{0} \vartriangleleft 
\Cal U ^{\alpha}_{1}
\vartriangleleft \cdots \vartriangleleft \Cal U^{\alpha}_{\beta} 
\vartriangleleft
\cdots \; \; \; (\beta<o^{\Cal U}(\alpha))$.
The following is proved in [Mi83].

\proclaim{Proposition 1.3} There is a class sequence $\overrightarrow{\Cal U}$
such that 
$L[\overrightarrow{\Cal U}] \models$ `` For every $\alpha$ $\overrightarrow
{\Cal U}\restriction (\alpha +1)$ is a coherent sequence,
every measure on $\alpha$ is equal to some $\Cal U ^{\alpha}_{\beta}$
and $o^{\Cal U}(\alpha) = \min \{ (\text{Mitchell order})^{V}(\alpha), \; 
\alpha^{++} \}$'' . Moreover,
$L[\overrightarrow{\Cal U}]$  satisfies GCH.
\endproclaim

We say that $\k$ has {\it a repeat point} (see [Mi82]) if there is a coherent 
sequence
$\overrightarrow{\Cal U}$ up to $\k$ and an ordinal $\theta < o^{\Cal U}(\k)$
such that
$$ \forall X \in \Cal U^{\k}_{\theta} \; \exists \alpha < \theta \; : \; X \in 
\Cal U^{\k}_{\alpha}\;.$$
It can be  proved that such $\theta$ must be greater than $\k^{+}$. Suppose we 
have a coherent
sequence $\overrightarrow{\Cal U}$ such that $o^{\Cal U}(\k)=\k^{++}$; then
using a simple counting argument we can prove the existence of a repeat
point for $\k$. Consequently, if Mitchell order of $(\k)$ is $\k^{++}$
then there is an inner model
satisfying GCH where $\k$ has a repeat point.

Our result is the following:

\proclaim {Theorem} If $\k$ has a repeat point in the ground model $V$ 
satisfying
GCH then there is  a generic extension of $V$ preserving cardinalities, 
cofinalities
and GCH in which Full Reflection holds at all $\lambda \leq \k$ and $\k$
is measurable.
\endproclaim

Actually if we start with $V=L[\overrightarrow{\Cal U}]$ from proposition 1.3
then our construction provides a class generic extension $V[G]$ preserving
cardinalities, cofinalities and GCH such that for any cardinal $\lambda$
$V[G]$ satisfies Full Reflection at $\lambda$, (Mitchell order)$^V (\lambda)
=o^{V[G]}(\lambda)$ and if $\lambda$ has a repeat in $V$ then $\lambda$ is
measurable in $V[G].$

\subhead
2. The Forcing $P_{\k+1}$
\endsubhead

From now on we work in a ground model $V$ satisfying GCH with a coherent
sequence $\overrightarrow{\Cal U}$ up to $\k$ and a repeat point at $\k$. For
$\lambda \leq \k$ let $\theta_{\lambda}$ be  $o^{\Cal U}(\lambda)$ 
if $\lambda$ does not have a repeat point, or otherwise the least $\theta$
such that $\Cal U^{\lambda}_{\theta}$ is a repeat point.

As usual, if $P$ is a forcing notion then $V(P)$ denotes either
the Boolean valued model or a generic extension by a $P$-generic 
filter over  $V$.

$P_{\k+1}$ will be an Easton support iteration of $\langle Q_{\lambda};\; 
\lambda \leq\k \rangle$, $Q_{\lambda}$
will be nontrivial only for $\lambda$ Mahlo.
$Q_{\lambda}$ (for $\lambda$ Mahlo) is defined in $V(P_{\lambda})$,
where $P_{\lambda}$ denotes the iteration below $\lambda$, 
as an iteration of length
$\lambda^{+}$ with $<\lambda$-support of forcing notions shooting clubs
through certain sets $X\subseteq \lambda$ (we will denote this standard forcing
notion $CU(X)$), always with the property that $X\supseteq Sing(\lambda)$.
This condition will guarantee $Q_{\lambda}$ to be essentially $<\lambda$-closed
(i.e. for any $\gamma < \lambda$ there is a dense $\gamma$-closed subset of 
$Q_{\lambda}$).
$Q_{\lambda}$ will also satisfy the $\lambda^{+}$-chain condition. Consequently
$P_{\lambda}$ will satisfy $\lambda$-c.c. and will have size $\lambda$. 
Cardinalities, cofinalities and GCH will be preserved, stationary subsets of 
$\lambda$
can be made nonstationary only by the forcing at $\lambda$, not below $\lambda$,
and not after the stage $\lambda$ - after stage $\lambda$ no subsets of 
$\lambda$ are added.

We use the $\lambda^{+}$-chain condition  of $Q_{\lambda}$ to get a canonical
enumeration of length $\lambda^+$ of all the $\lambda^{+}$ $Q_{\lambda}$-names 
for subsets of $\lambda$ so that the $\beta$th name appears in $V(P_{\lambda} 
\ast Q_{\lambda} | \beta)$. Moreover for $\delta < \theta_{\lambda}$ we will 
define filters $F_{\delta}=F^{\lambda}_{\delta}$ in $V(P_\lambda \ast Q_\lambda 
| \beta)$.
Their definition will not be absolute, however the filters will extend the
$V$-measures $\; \Cal U^{\lambda}_\delta$ and will 
increase coherently during the iteration.         

\proclaim{Definition} An iteration $Q$ of $\langle CU(B_{\alpha});\; 
\alpha < \alpha_0
\rangle$ with $<\lambda$-support and length $\alpha_0 < \lambda^{+}$ is called
{\it an iteration of order} $\delta_0$  if for all $\alpha < \alpha_0,$
$$ V(P_{\lambda} \ast Q|\alpha) \models \; B_\alpha \in F_\delta \text{ for any }
\delta < \delta_0 \text{ and } Sing(\lambda) \subseteq B_{\alpha}\;\;\; .$$
\endproclaim

(Note that an iteration of order $\delta_0$ is also an iteration of order
$\delta,$ for all $\delta<\delta_0.)$

$Q_{\lambda}$ is then defined as an iteration of
$\langle CU(B_\alpha);\; \alpha < \lambda^{+} \rangle$
with $<\lambda$-support and length $\lambda^{+}$ so that every
$Q_\lambda | \alpha$ is an iteration 
of order $\theta_\lambda$ and all potential names  $\dot X \subseteq \lambda$
are used cofinally many times in the iteration as some $B_\alpha$.

Observe that $Q_\lambda$ can be represented in
$V(P_\lambda)$ as a set of sequences of closed bounded subsets of
$\lambda$ in $V(P_\lambda)$ rather than in $V(P_\lambda \ast Q_\lambda 
| \alpha)$.
Moreover if $\dot q$ is a $P_\lambda$-name such that $1 \force_{P_\lambda}
\dot q \in Q_\lambda$ then using the $\lambda$-chain condition of $P_\lambda$
there is a set $A \subseteq \lambda^{+}$ (in $V$) of cardinality $<\lambda$
and $\gamma_0 < \lambda$ so that 
$$1\force_{P_\lambda} \;\; \text{supp}\;\dot q
\subseteq A  \text{ and } \forall \alpha \in A : \; \dot q (\alpha)
\subseteq \gamma_0 \;\;.$$
Consequently, $Q_\lambda$ can be represented as a set of functions
$g:A \times \gamma_0 \rightarrow [P_\lambda]^{<\lambda}$ where 
$A \subseteq \lambda^{+}$, $|A|<\lambda$ and $\gamma_0 < \lambda$.
In this sense $Q_\lambda$ has cardinality $\lambda^{+}$ and any
$Q_\lambda|\alpha$ has cardinality at most $\lambda$.

We will need to lift various elementary embeddings to generic extensions.
For a review of basic methods see [WoC92]. We will often use the
following simple fact: Let $N$ be a submodel of $M$ such that 
$M \cap \; ^{\k}N \subseteq N$ and let $G$ be a filter $P$-generic$/M$
 where $P$ satisfies $\k^{+}$-c.c.  Then $$M[G] \cap
\;^{\k}N[G] \subseteq N[G]\;\;.$$
Moreover if $Q$ is $\k$-closed and $H$ is $Q$-generic$/M[G]$
then 
$$M[G \ast H] \cap \;^{\k}N[G \ast H] \subseteq N[G \ast H ] \;\; .$$

\subhead
Definition of filters $F_\delta$
\endsubhead

The filters $F_\delta$ in $V(P_\lambda \ast Q)$, where $Q$ is any iteration of 
order $\delta+1$, are defined by induction so that the following is satisfied:

\proclaim{Proposition 2.1} Let $Q$, $Q'= Q\ast R$ be two iterations
of order $\delta'+1$ then
$$F_{\delta'}^{V(P_\lambda \ast Q)} \; = \; F_{\delta'}^{V(P_\lambda
\ast Q')} \; \cap \; V(P_\lambda \ast Q) \; \; .$$
Moreover  $\; F_{\delta '}^{V(P_\lambda)}\cap V = \Cal U ^{\lambda}_{\delta'}$.
\endproclaim

\proclaim{Proposition 2.2} Let $j=j_{\delta'}$ be the canonical embedding
from $V$ into $V^{\lambda}{/\Cal U ^{\lambda}_{\delta'}} = M$ and
$Q$  an iteration of order $\delta'+1$. Then $j$ can be lifted to an elementary 
embedding from the generic extension $V(P_\lambda \ast Q)$ of $V$ to
a generic extension $M(jP_\lambda \ast jQ)$ of $M$.
\endproclaim

\proclaim{Lemma 2.3} Let $ N=V^{\lambda}{/\Cal U ^{\lambda}_\beta}$ for some
$\beta > \delta'$ and $Q$ be an iteration of order $\delta'+1$. Then
$$F_{\delta'}^{V(P_\lambda \ast Q)} \; = \; F_{\delta'}^{N(P_\lambda \ast Q)}
\; \; . $$
Note that it also means that the definition of $F_{\delta'}$
relativized to $N(P_\lambda \ast Q )$ makes sense.
\endproclaim

\proclaim{Lemma 2.4} Let $j=j_{\delta'}:V \rightarrow M$. Then any iteration $Q$ 
of order
$\delta'$ is an subiteration of $(jP_\lambda)^{\lambda}$,
where $(jP_\lambda)^{\lambda}$ is the factor of
$jP_\lambda = P_\lambda \ast (jP_\lambda)^{\lambda} \ast 
(jP_\lambda)^{>\lambda}.$ Consequently
for any $G^{*}$ $jP_\lambda$-generic$/V$ and any $q \in Q$ there is an
$H \in M[G^{*}]$ $Q$-generic$/V[G]$ containing $q$
given by an embedding of $Q$ as a subiteration of $(jP_\lambda)^{\lambda}$, 
where $G=G^{*} \restriction
P_\lambda$.
\endproclaim

\proclaim{Definition}
Let $j,\; Q,\; G^{*},\; G$ be as in the lemma. Then
$Gen_j (Q,G^{*})$ is the set of all filters $H\in M[G^{*}]$ 
$Q$-generic$/V[G]$ given by an embedding of $Q$ as 
a subiteration of $(jP_\lambda)^{\lambda}.$
\endproclaim

\proclaim{Lemma 2.5} Let $j$ be as above, $Q$ an iteration of order
$\delta'+1$, $G^{*}$ $jP_\lambda$-generic$/V$, $H \in Gen_j (Q,G^{*})$.
For every $\beta < l(Q)$ let $C_\beta \subset \lambda$ be the club
$\bigcup\{r_\beta; r\in H\}$, and let
$[H]^{j}$ denote the $j(l(Q))$-sequence given by
$$[H]^{j}_\gamma = \cases C_\lambda\cup \{ \lambda \} , & 
\text{ if } \gamma=j(\beta) \\
                           \emptyset & \text{ otherwise.} \endcases
$$
Then $[H]^{j} \in jQ_{/G^{*}}$.
\endproclaim

Propositions 2.1 and 2.2 will be essential to prove Full Reflection
in the generic extension. We will later prove that if lemma 2.3 holds
for $\delta' < \delta$ then lemma 2.4 holds for all $\delta' \leq \delta$.

Now suppose that the filters $F^{\lambda'}_{\delta'}$ 
were defined for all $\lambda'< \lambda$ and $\delta' < \theta_{\lambda'}$
and for $\lambda'=\lambda$ and $\delta'<\delta$ so that
2.1-2.5 holds. Moreover
let $\alpha < \lambda^{+}$ and $F_\delta=F^{\lambda}_{\delta}$ be defined for 
all iterations
of order $\delta+1$ of length $<\alpha$ so that lemma 2.5 holds for $\delta$
and iterations of length $\leq \alpha$. Let $j=j_\delta:V \rightarrow
M=V^{\lambda}{/\Cal U ^{\lambda}_\delta}$. Then we can define $F_\delta$
for iterations of order $\delta+1$ and length $\alpha$.

\proclaim{Definition} Let Q be an iteration of order $\delta+1$ and length
$\alpha$, $j=j_\delta :V \rightarrow M.$ For a $P_\lambda \ast Q$-name 
$\dot X$ of a subset of $\lambda$
and $(p,q) \in P_\lambda \ast Q$ we define
$$ (p,q) \force _{P_\lambda \ast Q} \; \dot X \in F_\delta $$
if the following holds in V:

$$j(p) \force_{jP_\lambda}\text{ ``For any $H \in Gen_j (Q,G^{*})$ containing 
$q$ :} $$
$$[H]^{j} \force_{jQ} \; \check \lambda \in j\dot X \;\;\; \text{''.}$$
\endproclaim

The definition says that $(p,q) \force \dot X \in F_\delta$ if 
$\lambda \in j^{*} X$ whenever $j^{*}:V[G\ast H] \rightarrow M[G^{*} 
\ast H^{*}]$
is a lifting of $j$ of certain kind and $(p,q)\in G\ast H$. To verify
soundness of the definition let us first prove lemma 2.4 for $\delta$.

\demo{Proof of lemma 2.4} 
Let $j=j_\delta : V \rightarrow M=V^{\lambda}{/\Cal U ^{\delta}_{\lambda}}$,
Q be an iteration of order $\delta$. We assume that lemma 2.3 holds for
$\delta'<\delta$. Observe that $jP_{\lambda} = P_{\lambda} \ast
(jP_\lambda)^{\lambda} \ast (jP_\lambda)^{>\lambda}$ is an Easton support
iteration (in $M$) below $j\lambda$ and $(jP_\lambda)^{\lambda}$ is an iteration
of length $\lambda^{+}$ with $<\lambda$-support
such that for any $\alpha < \lambda^{+}$ $(jP_\lambda)^{\lambda}|\alpha$
is an iteration of order $\delta = \theta_\lambda^{M}$ and all potential
names $\dot X \subseteq \lambda$ are used cofinally many times in the
iteration. That is true in $M(P_\lambda)$ as in $V(P_\lambda)$.

Let us now define what it means for $Q$ to be a subiteration of
$(jP_\lambda)^\lambda .$  Suppose that $P$,
$Q$ are iterations of legths $l(P) \leq l(Q)$ of 
$\langle \dot R_\gamma; \; \gamma<l(P) \rangle$ and
$\langle \dot S_\alpha; \; \alpha<l(Q) \rangle$ with $<\lambda$-support
essentially $<\lambda$-closed. Then we say that $P$ is {\it a subiteration of}
$Q$ if there is an increasing sequence 
$\langle \alpha_\gamma;\; \gamma<l(P) \rangle$ of ordinals below
$l(Q)$ such that
$$\align
1) \;\;\; &Q_{\alpha_0} \force \check R_0 = \dot S_{\alpha_0}\\ 
&\text{consequently for }\beta> \alpha_0 \;\; Q_\beta \simeq
P_1 \ast Q_\beta '\;,\\
2) \;\;\; &\text{and by induction }Q_{\alpha_\gamma}\force
\check {\dot R} _{\gamma} = \dot S_{\alpha_\gamma}\;,\\
&\text{using the inductive assumption that }Q_{\alpha_\gamma} \simeq
P_\gamma \ast Q_{\alpha_\gamma}'\;,\\
&\text{consequently again for }\beta > \alpha_\gamma\text{ we have }
Q_\beta \simeq P_{\gamma+1}\ast Q_\beta ' \; \; .
\endalign
$$
It is now obvious that in this sense $Q$ is a subiteration of
$(jP_\lambda)^{\lambda}$. Moreover for any $\alpha<\lambda^{+}$ there
is a sequence $\langle \alpha_\gamma;\; \gamma < l(Q) \rangle$
determining an embedding of $Q$ into $(jP_\lambda)^{\lambda}$ such that
$\alpha_0 > \alpha$.

Finally let $G^{*}$ be $jP_\lambda$-generic$/V$ and $q\in Q$.
Then $(G^{*})^{\lambda}$ is $(jP_\lambda)^{\lambda}$-generic$/V[G]$.
Note that the set
$$\align
D=\{ r\in (jP_\lambda)^{\lambda};&\text{ there is a sequence }
\langle \alpha_\gamma; \gamma<l(Q) \rangle \\
&\text{ determining an embedding of }
Q \text{ into } (jP_\lambda)^{\lambda} \\ 
&\text { such that } q 
\text{ corresponds to } r\restriction \langle \alpha_\gamma;
\gamma<l(Q) \rangle \;\;\;\}
\endalign
$$
is dense in $(jP_\lambda)^{\lambda}$. Thus let $r\in D \cap (G^{*})^{\lambda}$
and $\langle \alpha_\gamma;
\gamma<l(Q) \rangle $ be the sequence. Then $(G^{*})^{\lambda} \restriction 
\langle \alpha_\gamma;
\gamma<l(Q) \rangle $ gives the $Q$-generic$/V[G]$ filter $H\ni q$.
\qed
\enddemo

We have defined $F_\delta$ for iterations of length $\leq \alpha$.
Let us now prove lemma 2.5 for iterations of length $\alpha+1$.

\demo{Proof of lemma 2.5} 
Let $j=j_\delta:V\rightarrow M$, $Q$ be an iteration of order $\delta+1$
of $\langle CU(B_\beta); \beta< \alpha+1 \rangle,$
 $G^{*}$ $jP_\lambda$-generic$/V$, $H\in Gen_j (Q,G^{*})$,
 $[H]^{j}$ as in the lemma. Note that $[H]^{j}$
is a $j(\alpha+1)$-sequence of closed bounded subsets of $j\lambda$,
$\text{supp}\;[H]^{j}=j"(\alpha +1)$. 
Since $|j"(\alpha +1)|=\lambda < j\lambda$
in $M$, $[H]^{j}$ has a small support and $[H]^{j} \in M[G^{*}].$ It follows 
from the induction
hypothesis that $[H]^{j} \restriction j(\alpha) \in j(Q\restriction \alpha)$.
We only have to verify that in $M[G^{*}]$
$$
[H]^{j}\restriction j(\alpha) \force_{j(Q)\restriction j(\alpha)}
\; [H]^{j}_{j(\alpha)} \in j(\langle CU(B_\beta); \beta< \alpha+1 
\rangle)_{j(\alpha)}
$$
which means just that
$$
[H]^{j}\restriction j(\alpha) \force [H]^{j}_{j(\alpha)}\; \subseteq \;
j(B_\alpha)
$$
where $[H]^{j}_{j(\alpha)}= \cup\{r_\alpha;\;r\in H\}\cup\{\lambda \}$.
If $\gamma \in [H]^{j}_{j(\alpha)}$, $\gamma < \lambda$, then $\gamma\in
r_\alpha$ for some $r\in H$ and 
$$V[G] \models r\restriction \alpha \force_{Q|\alpha}
r_\alpha \subseteq B_\alpha$$
so
$$ M[G^{*}] \models j(r\restriction \alpha) \force_{j(Q|\alpha)} 
r_{\alpha}=jr_\alpha
\subseteq jB_\alpha \;\; . $$
Since $[H]^{j}\restriction j(\alpha)\; \leq \; j(r\restriction \alpha)$
we get $[H]^{j}\restriction j(\alpha) \force \gamma \in jB_\alpha\; .$
So we only have to prove that in $M[G^{*}]$
$$ [H]^{j}\restriction j(\alpha) \force \lambda \in jB_\alpha \tag*$$
Here we use the fact that $Q$ is an iteration of order $\delta+1$.
It implies that $Q\restriction\alpha \force B_\alpha \in F_\delta$
and that exactly gives (*) in $V[G^{*}]$ and so in $M[G^{*}]$ by the 
definition of $F_\delta$ in
$V(P_\lambda \ast Q|\alpha).$
\qed
\enddemo

$F_\delta$ is now well defined in $V(P_\lambda \ast Q)$ for any iteration
$Q$ of order $\delta +1$. We have to verify Proposition 2.1, 2.2 and Lemma
2.3 for $F_\delta.$

\demo{Proof of Lemma 2.3} 

Let $\beta>\delta$, $N=V^{\lambda}{/\Cal U^{\lambda}_{\beta}}$, 
$j=j_\delta : V \rightarrow M$, $Q$ an iteration of order $\delta +1.$
We want to prove that
$$ F_\delta^{V(P_\lambda\ast Q)} = F_\delta^{N(P_\lambda \ast Q)} \;\;.$$
Let $j':N\rightarrow N^{\lambda}{/\Cal U^{\lambda}_\delta} = N'$, observe
that $j'=j\restriction N$. We need to prove that the following two conditions
are equivalent:
$$\align
V \models jp \force_{jP_\lambda} &\forall H \in Gen_j (Q,G^{*}),
\;H\ni q \tag 1\\
&[H]^{j} \force _{jQ} \check \lambda \in j\dot X \;\;,\\
N \models  jp \force_{jP_\lambda} &\forall H \in  Gen_j (Q,G^{*}),
\;H\ni q \tag 2\\
&[H]^{j} \force _{jQ} \check \lambda \in j\dot X\;\; .
\endalign
$$
Observe that $V[G^{*}]\cap\;^{\lambda}N[G^{*}]\;\subseteq \; N[G^{*}],$
and so $$Gen_j (Q,G^{*})^{V[G^*]} = Gen_j (Q,G^{*})^{N[G^*]} .$$   From that
the equivalence of (1) and (2) easily follows.
\qed
\enddemo

\demo{Proof of Proposition 2.2} Let $j=j_\delta : V \rightarrow M,$
$Q$ be an iteration of order $\delta +1$. Let $G^{*}$ be 
$jP_\lambda$-generic$/V,$
$G=G^{*}\restriction P_\lambda.$ Then $j$ is lifted to
$j^{*}:V[G] \rightarrow M[G^{*}]\;.$
Using Lemma 2.4 find $H\in M[G^{*}]$ $Q$-generic$/V[G].$ By Lemma 2.5
$[H]^{j}\in j^{*}Q$ and $r\in H$ implies $j^{*}r\geq [H]^{j}.$ Thus let
$H^{*}$ be $jQ$-generic$/V[G^{*}]$ containing $[H]^{j}.$ Then $r\in H$ implies
$j^{*}r \in H^{*}$ and $j^{*}$ is lifted to $j^{**}:V[G\ast H] \rightarrow
M[G^{*} \ast H^{*}]\; .$
\qed
\enddemo

\demo{Proof of Proposition 2.1}
Let $Q$, $Q'=Q\ast R$ be iterations of order $\delta +1$, we want to prove
$$F_{\delta}^{V(P_\lambda \ast Q)} \; = \; F_{\delta}^{V(P_\lambda
\ast Q')} \; \cap \; V(P_\lambda \ast Q) \; \; .$$

Firstly let us prove the easy direction: $\dot X \subseteq \lambda$
a $P_\lambda \ast Q$-name, $(p,q)\in P_\lambda\ast Q,$ 
$(p,q) \force_{P_{\lambda}\ast Q} \dot X \in F_\delta\; .$
Then it is straightforward  that $(p,q^{\frown}1)\force_{P_{\lambda}\ast Q'}\;
\dot X \in F_\delta\; .$

Now let $\dot X \subseteq \lambda$ be a $P_\lambda \ast Q$-name, 
$(p,q')\in P_\lambda\ast Q'$ and $(p,q') \force_{P_{\lambda}\ast Q'} 
\dot X \in F_\delta\; .$
We prove that $(p,q) \force_{P_{\lambda}\ast Q} \dot X \in F_\delta$ where
$q=q' \restriction l(Q).$
 Let $G^{*}\ni p$ be $jP_\lambda$-generic$/V,$ 
$j=j_\delta : V \rightarrow M,$ $H\in  Gen_j (Q,G^{*})$
and $q\in H.$ We want to prove that $[H]^{j}\force_{j(Q)} \check \lambda
\in j\dot X \; .$ Suppose not, then there is $\tilde q\leq [H]^{j}$
such that $\tilde q\force_{j(Q)} \check \lambda \notin j\dot X\; .$
Obviously there exists $\tilde H \in  Gen_j (Q',G^{*})$ such
that $\tilde H \restriction Q = H$ and $ \tilde H \ni q' .$
Note that $\tilde q^{\frown}1$ and $[\tilde H]^{j}$ are compatible,
 $\tilde q^{\frown}1 \cup
[\tilde H]^{j} \in jQ'.$
But $[\tilde H]^{j} \force _{jQ'} \check \lambda \in j\dot X$
and $\tilde q^{\frown}1 \force_{jQ'} 
\check \lambda \notin j\dot X$ - a contradiction.

Finally let us prove that $\Cal U ^{\lambda}_{\delta} = F^{V(P_\lambda)}_\delta
\cap V .$
The inclusion $\Cal U^{\lambda}_{\delta} \subseteq F_\delta$ is obvious.
Now let $X \in V,$ $X\subseteq \lambda$ and $p\force _{P_\lambda} 
\check X \in F_\delta.$ Then it means that $j(p) \force _{jP_\lambda} 
\check \lambda \in
j(\check X)$ which can be true only if $\lambda \in j(X).$ So
$X \in \Cal U^{\lambda}_{\delta}.$ 
\qed
\enddemo

\subhead
3. Full Reflection in $V(P_{\k +1})$
\endsubhead

To prove that Full Reflection holds in
$V(P_{\k +1})$ at some $\lambda \leq \k$ it is 
enough to prove that in $V(P_{\lambda +1})$.
Fix $\lambda \leq \k.$

Firstly let us prove the existence of sets $E_\delta$ $(\delta<\theta_\lambda)$
separating the measures $\Cal U^{\lambda}_\delta$ in the sense that
$E_\delta \in \Cal U^{\lambda}_{\delta'}$ iff $\delta=\delta' .$

\proclaim{Proposition 3.1}
Suppose that $\Cal U^{\lambda}_\delta$ is not a repeat point. Then
there is a set $E \in \Cal U^{\lambda}_\delta$ such that $E \notin 
\Cal U^{\lambda}_{\delta'}$ for any $\delta' \ne \delta,$ $\delta' < 
o^{\Cal U}(\lambda).$
\endproclaim
\demo{Proof}
Since $\Cal U^{\lambda}_\delta$ is not a repeat point there is a set
$X \in \Cal U^{\lambda}_\delta$ such that $X \notin \Cal U^{\lambda}_{\delta'}$
for all $\delta' < \delta.$ Define for $\xi < \lambda$
$$
f_\delta(\xi )=\sup\{\eta;\;\eta \leq o^{\Cal U}(\xi) \;\and \; \forall
\eta'<\eta:X\cap\xi\notin\Cal U^{\xi}_{\eta'} \}
$$
and
$$
Y=\{ \xi<\lambda;\;f_\delta(\xi)=o^{\Cal U}(\xi) \}.
$$
\proclaim{Claim} $Y \in \Cal U^{\lambda}_\delta$ but $Y\notin \Cal 
U^{\lambda}_{\delta'}$ for any $\delta'$ such that $\delta<\delta'<
o^{\Cal U}(\lambda).$
\endproclaim
\demo{Proof}
1) Let $j:V \rightarrow N=V^{\lambda}{/\Cal U^{\lambda}_\delta}.$
Then
$$
Y\in\Cal U^{\lambda}_\delta \text{ iff } \lambda \in jY = \{ \xi<j\lambda;\;
N\models o^{j\Cal U}(\xi)=j(f_\delta)(\xi) \}$$
$$\text{iff }N\models o^{j\Cal U}(\lambda) = j(f_\delta)(\lambda)$$
By the definition of a coherent sequence $o^{j\Cal U}(\lambda)=\delta.$
Moreover
$$
N \models j(f_\delta)(\lambda) = sup \{ \eta;\; \eta \leq \delta \and
\forall\eta'<\eta:\;jX\cap\lambda \notin j\Cal U^{\lambda}_{\eta'} \}.
$$
By  coherence $j\Cal U^{\lambda}_{\eta'} =\Cal U^{\lambda}_{\eta'}$
for $\eta'<\delta$, and also $jX\cap \lambda=X.$
Consequently $j(f_\delta)(\lambda)=\delta$ and $Y\in \Cal U^{\lambda}_{\delta}.$

2) Let $\delta<\delta'<o^{\Cal U}(\lambda)$ and $j:V \rightarrow N=
V^{\lambda}{/\Cal U^{\lambda}_{\delta'}}.$
Again $Y\in \Cal U^{\lambda}_{\delta'}$ iff $N\models o^{j\Cal U}(\lambda)=
j(f_\delta)(\lambda).$
As above $o^{j\Cal U}(\lambda)=\delta'$ but $j(f_\delta)(\lambda)=\delta$
since $X\notin \Cal U^{\lambda}_{\delta}.$
Thus $Y\notin \Cal U^{\lambda}_{\delta'}$ and the claim is proved.
\enddemo
Finally put $E=Y\cap X.$
\qed
\enddemo
So we have a separating sequence. We can suppose that the $E_\alpha$ are sets
of regular cardinals because $Reg(\lambda)\in \Cal U^{\lambda}_{\delta}$
for any $\delta.$ We are going to prove the following:

\proclaim{Proposition 3.2}
In $V(P_{\lambda +1})$ the sets $\langle E_\alpha;\;\alpha<\theta_\lambda 
\rangle$
form a maximal antichain of stationary subsets of $Reg(\lambda).$
Moreover if $S\subseteq E_{-1}=Sing(\lambda)$ or $S\subseteq E_\alpha$
is stationary then
S reflects in any $E_\beta$ for $\beta>\alpha$ 
and $Tr(E_\alpha)\cap E_\alpha$ is nonstationary if $\alpha>-1$.
Consequently each $E_\alpha$ is a canonical stationary set
of order $\alpha$, $o(\lambda)=\theta_\lambda$ and Full Reflection at
$\lambda$ holds.
\endproclaim

To prove the proposition we need the following lemmas:

\proclaim{Lemma 3.3}
$V(P_\lambda \ast Q_\lambda |\alpha) \models Club(\lambda) \subseteq F_\delta$
for any $\delta < \theta_\lambda . $
\endproclaim

\demo{Proof}
Let $\delta < \theta_\lambda$ and $(p,q) \force_ {P_\lambda \ast 
Q_\lambda |\alpha} \; \text{`` }\dot X \subseteq \lambda \text{ is a club''} .$
If $(p,q)$ does not force $\dot X \in F_\delta$ then
by the definition of $F_\delta$ there is a $G^{*}$ $jP_\lambda$-generic$/V,$
$G^{*}\ni jp,$ where $j:V \rightarrow V^{\lambda}{/\Cal U^{\lambda}_\delta}
=M,$ and $H \in  Gen_j (Q_\lambda|\alpha,G^{*}),$ $H\ni q,$
such that in $V[G^{*}]$ $[H]^{j} \notforce_{j(Q_\lambda|\alpha)}\; 
\check \lambda \in j\dot X.$ 
So there is $H^{*} \ni [H]^{j}$ $j(Q_\lambda |\alpha)$-generic
$/V[G^{*}]$ so that $\lambda \notin j\dot X_{/G^{*} \ast H^{*}}.$ 
The embeding $j$ is by the proof of Propostion 2.2 lifted to
the elementary embeding $j^{**}:V[G\ast H] \rightarrow M[G^{*} \ast H^{*}],$
$X=\dot X_{/G\ast H}$ is a club in $V[G\ast H]$, thus $j^{**}X=j\dot X_{/G^{*} 
\ast H^{*}}$ is a club in $M[G^
{*} \ast H^{*}].$ Since $j^{**}X\cap \lambda =X$
necessarily $\lambda \in j^{**}X$ - a contradiction.
\qed
\enddemo

\proclaim{Lemma 3.4} 
If $S$ in $V(P_\lambda \ast Q_\lambda |\alpha)$
is $F_\delta$-positive then $Tr(S) \in F_{\delta'}$ for any
$\delta' > \delta .$ 
If $S \subseteq Sing(\lambda)$ is stationary then $Tr(S) \in F_\delta$ for 
any $\delta .$
Moreover $V(P_\lambda) \models E_\delta \setminus 
Tr(E_\delta) \; \in \; F_\delta.$
\endproclaim

\demo{Proof}
Suppose $(p,q) \force _{P_\lambda\ast Q_\lambda |\alpha} \; \text{``$\dot S$ 
is $F_\delta$-positive''}$ but
$(p,q) \notforce _{P_\lambda\ast Q_\lambda |\alpha} \;$``$Tr(\dot S) \in 
F_{\delta'}$'' for some $\delta' > \delta.$
As usual denote $j:V \rightarrow M=V^{\lambda}{/\Cal U^{\lambda}_{\delta'}} .$
Then there is a filter $G^{*}$ $jP_\lambda$-generic$/V,$ $G^{*} \ni jp,$
and $H \in  Gen_j (Q_\lambda|\alpha,G^{*}),$ $H \ni q,$
and a filter $H^{*}$ $j(Q_\lambda |\alpha)$-generic$/V[G^{*}],$
$H^{*} \ni [H]^{j} ,$
so that $\lambda \notin jTr(\dot S)_{/G^{*} \ast H^{*}} \; .$
As above $j$ is lifted to
$$ j^{**}: V[G\ast H] \rightarrow M[G^{*} \ast H^{*}] \; . $$
$S=\dot S_{/G\ast H}$ is $F_\delta$-positive in $V[G\ast H]$ and
$$jTr(\dot S)_{/G^{*} \ast H^{*}} \; = \; j^{**}(Tr(S)) \; = \; 
Tr^{M[G^{*}\ast H^{*}]}(j^{**}S)\; .$$
Thus $\lambda \notin Tr^{M[G^{*}\ast H^{*}]}(j^{**}S)$ which means that
$$ M[G^{*} \ast H^{*}] \models \text{`` $S$ is not stationary in $\lambda$''}$$
because $S=j^{**}S \cap \lambda.$ Consequently also
$$V[G^{*}\ast H^{*}] \models  \text{``$S$ is not stationary in $\lambda$''}.$$
Observe that $j(Q_\lambda |\alpha)$ has a dense subset $\lambda$-closed
in $M[G^{*}]$ and thus also in $V[G^{*}].$ Moreover
$jP_\lambda = P_\lambda \ast (jP_\lambda)^{\lambda} \ast R$ where $R$
is essentially $\lambda$-closed in $V[G^{*} | \lambda +1].$ It implies that
already $V[G^{*} | \lambda +1] \models  \text{``$S$ is not stationary 
in $\lambda$''}.$
Let us now consider the isomorphism $(jP_\lambda)^{\lambda} \simeq
(Q_\lambda |\alpha) \ast \tilde Q$ from the proof of Lemma 2.4 giving
the filter $H=G^{*} \restriction (Q_\lambda |\alpha),$
let $\tilde H = G^{*} \restriction \tilde Q.$ Since every subset of $\lambda$
in $V[G\ast H \ast \tilde H]$ is already in some $V[G\ast H \ast 
\tilde H|\beta]$ there is a $\beta < \lambda^{+}$ so that
$$V[G\ast H \ast \tilde H|\beta] \models \text{``$S$ is not stationary 
in $\lambda$''}.$$
But since $(Q_\lambda |\alpha) \ast (\tilde Q |\beta)$ is an iteration of 
order $\delta' \geq \delta+1$
it follows from Proposition 2.1 that
$$V[G\ast H \ast \tilde H|\beta] \models \text{``$S$ is $F_\delta$-positive''}$$
which contradicts Lemma 3.3.

The proof for $S \subseteq Sing(\lambda)$ is the same using the following
fact instead of Proposition 2.1.

\proclaim{Claim}
Stationary subsets of $Sing(\lambda)$ are preserved by iterations of
order 0.
\endproclaim
\demo{Proof}
For simplicity assume that $R=CU(X)$ where $X\supseteq Sing(\lambda) ;$
the generalization for an iteration of order 0 is straightforward. We closely
follow the proof of 7.38 in [J86].

Let $S\subseteq Sing(\lambda)$ be stationary, $\dot C$ an $R$-name and
$p\force_R\;\dot C \subseteq \lambda$ is a club. We need a $\tilde q \leq p$
and $\beta \in S$ so that $\tilde q\force \beta \in \dot C .$ Put $A_0 = 
\{ p\},$
$\gamma_0 = \max(p),$ and inductively for $q\in A_\alpha$ find $r(q) \leq q$
and $\beta(q) > \gamma_\alpha$ so that $\max(r(q)) > \gamma_\alpha$ and
$r(q)\force \beta(q) \in \dot C .$ Put
$$ A_{\alpha +1} = A_\alpha \cup \{r(q);\;q\in A_\alpha\} \text{ and}$$
$$\gamma_{\alpha +1} = \sup (\{\max(q);\; q\in A_{\alpha+1}\} 
\cup \{\beta(q);\;
q\in A_\alpha \}).$$
For $\beta$ limit put
$$\align
A_\beta =\bigcup_{\alpha < \beta} A_\alpha \;\cup\;\{ &\text{unions of all 
decreasing sequences}\\
& \subseteq \bigcup_{\alpha < \beta} A_\alpha 
\text{ that are in } R \}\text{ and}
\endalign
$$
$$\gamma_\beta = \sup \{\gamma_\alpha ;\; \alpha < \beta \} .$$
Find a $\beta \in S$ such that $\gamma_\beta = \beta.$ Observe that 
cf$(\beta)<\beta$ and all unions of increasing sequences $\subseteq 
\bigcup_{\text{cf}(\beta) < 
\alpha < \beta} A_\alpha$ of length $\leq cf(\beta)$
are in $R.$ Now it is easy to find an increasing sequences $\beta_\alpha
\nearrow \beta$ and decreasing $q_\alpha \searrow \tilde q \in R$
($\alpha < cf(\beta)$) so that $q_\alpha \force \beta_\alpha \in \dot C .$
Consequently $\tilde q \force \beta \in \dot C.$
\enddemo
 
Let us now prove that $V(P_\lambda) \models E_\delta \setminus Tr(E_\delta) 
\in  F_\delta \; .$
Let $j=j_\delta:V \rightarrow M$ then $(jP_\lambda)^{\lambda}$ is an iteration
of length $\lambda^{+}$ such that $(jP_\lambda)^{\lambda} |\alpha$
is always an iteration of order $\delta$ and every potential name 
is used cofinally many times. Thus a club is shot through $\lambda \setminus 
E_\delta$ in the iteration. It implies that
$$V[G^{*}] \models E_\delta \subseteq \lambda \text{ is nonstationary}$$
and consequently
$$V[G^{*}] \models \lambda \in j(E_\delta \setminus Tr(E_\delta)).$$
\qed
\enddemo

\demo{Proof of Proposition 3.2}
That each $E_\delta$ is stationary in $V(P_{\lambda+1})$ folows easily
from Proposition 2.1 and Lemma 3.3. Let $\delta \ne \delta' < \theta_\lambda,$ 
then $\lambda \setminus (E_\delta \cap E_{\delta'}) \; \in \;
\Cal U^{\lambda}_\eta$ for any $\eta < \theta_\lambda$ and $\lambda \setminus 
(E_\delta \cap E_{\delta'}) \supseteq Sing(\lambda),$ so
$\lambda \setminus (E_\delta \cap E_{\delta'})$ contains a club in
$V(P_{\lambda +1}),$ and $E_\delta \cap E_{\delta'}$ is nonstationary.
Let now $A\subseteq Reg(\lambda),$ $A\in V(P_{\lambda +1})$ be such that
$A \cap E_\delta$ is nonstationary in $V(P_{\lambda +1})$ for any
$\delta < \theta_\lambda.$
We know that $A \in V(P_\lambda \ast Q_\lambda |\beta)$ for some $\beta < 
\lambda^{+}.$
\proclaim{Claim}
$V(P_\lambda \ast Q_\lambda|\beta) \models \lambda \setminus A \; \in \; 
F_\delta \text{ for any } \delta < \theta_\lambda .$
\endproclaim
\demo{Proof}
If $A$ was $F_\delta$-positive then $E_\delta \cap A$ would be 
$F_\delta$-positive in $V(P_\lambda \ast Q_\lambda|\alpha)$ for
$\alpha \geq \beta.$ Therefore $E_\delta \cap A$ would be stationary in 
$V(P_{\lambda+1}).$
\enddemo

Since also $Sing \subseteq \lambda \setminus A$ there is a 
club $C \subseteq \lambda 
\setminus A $ in $V(P_{\lambda +1})$, and so $A$ is nonstationary.

We have proved that $\langle E_\delta;\;\delta < \theta_\lambda \rangle$
forms a maximal antichain of stationary subsets of $Reg(\lambda)$ in
$V(P_{\lambda +1}).$

Now let $S\subseteq E_\delta$ be stationary, $\delta' > \delta,$
$S\in V(P_\lambda \ast Q_\lambda |\alpha).$ $S$ is $F_\delta$-positive
(or just stationary if $\delta=-1$)
in $V(P_\lambda \ast Q_\lambda |\alpha)$ and so by Lemma 3.4
$Tr(S) \in F_{\delta'}.$
Consequently $E_{\delta'} \setminus Tr(S)$ is nonstationary in
$V(P_{\lambda+1})$ - $(\lambda \setminus E_{\delta'})  \cup  (E_{\delta'}\cap 
Tr(S))$ contains a club - which exactly means that $S<E_{\delta'}.$
Since by Lemma 3.4 $E_\delta \setminus Tr(E_\delta) \; \in \; F_\delta$
the set $Tr(E_\delta)\cap E_\delta$ is nonstationary in $V(P_{\lambda+1})$
- $(\lambda \setminus E_\delta ) \cup ( E_\delta\setminus Tr(E_\delta))$
contains a club.
\qed
\enddemo

The following easy observation tells us more about the properties of the
algebra $\Cal P (\kappa) / NS $ in the resulting model.

\proclaim{Proposition 3.5}
Let $-1 \leq \alpha < \beta < \theta_\lambda ,$ then the sum of the sets
$\{ E_\delta ; \alpha < \delta \leq \beta \} $ in the algebra
$\Cal P (\kappa) / NS $ exists. Moreover for any normal measure over
$\kappa$ this sum has measure zero.
\endproclaim
\demo{Proof}
It follows immediately from proposition 3.2 that $Tr(E_\alpha)$ is
the sum of $\{ E_\delta ; \alpha < \delta < \theta_\lambda \} .$
Hence the desired sum is just $Tr(E_\alpha) \setminus Tr(E_\beta) .$
For any normal measure over $\kappa$ the measure of $Tr(E_\beta)$ is one,
consequently the measure of $Tr(E_\alpha)\setminus Tr(E_\beta)$ must be zero.
\qed
\enddemo

\subhead
4. Measurability of $\kappa$ in $V(P_{\k+1})$ 
\endsubhead

Let $\Cal U^{\k}_\theta$ be the first repeat point of $\k$ and
$j=j_\theta:V\rightarrow M = V^{\k}{/\Cal U^{\k}_\theta}.$
Then $(j\overrightarrow{\Cal U})\restriction \k+1\; =\; \overrightarrow{\Cal U}
\restriction (\k, \theta)$ and it follows from Lemma 2.3 that $(jP_\k)^{\k}$
is an iteration of length $\k^{+}$ with $<\k$-support such that any initial
segment is an iteration of order $\theta$ and any potential name is used
cofinally many times in $M(P_\k)$ as well as in $V(P_\k).$
Consequently we can suppose that $Q_\k = (jP_\k)^{\k}.$

Using methods for extending elementary embeddings
(see [WoC92] and [JWo85]) we will prove that $\kappa$ is actually
measurable in $V(P_{\k +1}).$ 
Let $G$ be a $P_\k$-generic filter$/V$, $G_\k$ a $Q_\k$-generic$/V[G].$
We know that $jP_\k = P_\k \ast Q_\k \ast R,$ 
where the factor $R$ is $\k$-closed
in $M[G][G_\k]$ and consequently in $V[G][G_\k].$

\proclaim{Lemma 4.1}
There is a filter $H\in V[G][G_\k]$ $R$-generic$/M[G][G_\k].$
\endproclaim
\demo{Proof}
Since $V\models |P_\k|=\k$ we have $M \models |jP_\k|=j\k$
and therefore the factor $R$ has cardinality $j\k$ in $M[G][G_\k].$
Thus $M[G][G_\k] \models |\Cal P(R)|=(j\k)^{+}$ because of GCH.
Put
$$ \Cal D = \{ D\in M[G][G_\k];\; D \text{ is a dense subset of } R \} $$
then the cardinality of $\Cal D$ in $V[G][G_\k]$ is same as the cardinality
of $(j\k)^{+M}$ which is $\k^{+}.$ Now use the fact that $R$ is $\k$-closed
to get a generic filter $H\in V[G][G_\k].$
\qed
\enddemo

Consequently $j$ can be lifted to
$$ j^{*}: V[G] \rightarrow M[G][G_\k][H]$$
where $j^{*}$ is defined in $V[G][G_\k].$
 Next we need to prove the following important lemma:

\proclaim{Lemma 4.2}
For any $\alpha < \k^{+}$
$$[G_\k \restriction \alpha]^{j} \; \in \; j^{*}(Q_\k |\alpha) $$
where $[G_\k \restriction \alpha]^{j}$ is defined as in Lemma 2.5.
\endproclaim
\demo{Proof}
Let $j_\delta$ denote the elementary embedding $j_\delta:V \rightarrow
M_\delta=V^{\k}{/\Cal U^{\k}_\delta}$ for $\delta < \theta.$
It follows from Lemma 2.5 and the proof of Lemma 2.3 that
for any $\delta<\theta$
$$\align
M_\delta \models 1\force_{j_\delta P_\k} \; & \forall H \in M_\delta[G^{*}]\;\;
Q_\k|\alpha\text{-generic}/M_{\delta}[G] \tag*\\
& [H]^{j} \in j_\delta(Q_\k|\alpha)_{/G^{*}}
\endalign
$$
Denote this formula $\varphi (j_\delta P_\k , Q_\k|\alpha , 
j_\delta(Q_\k|\alpha)).$

Now we need to introduce the notion of a canonical name. We say that
$f\in V^{\k}$ is a canonical name for $x \in V$ iff for any measure $\Cal U$
over $\k$ the set $x$ belongs to the transitive collapse $V^{\k}{/\Cal U}$
and is equal to $[f]_{\Cal U}.$ Let
$$ C=\{ x\in V ;\; x\text{ has a canonical name} \}.$$
Obviously $V_\k \subseteq C$ and $C^{\leq \k} \subseteq C.$ Since $P_\alpha \in 
V_\k$ for $\alpha < \k$ we get that $P_\k \in C$ and $Q_\k|\alpha \in C.$ Let f 
be the
canonical name for $Q_\k|\alpha$. Then by the Lo\'s 
Theorem (*) is equivalent to
$$\{\beta < \k;\; V\models \varphi(P_\k, f(\beta), Q_\k|\alpha) \} \in 
\Cal U^{\k}_\delta .$$
Since this is true for any $\delta < \theta$ and $\Cal U^{\k}_\theta$
is a repeat point it follows
$$\{\beta < \k;\; V\models \varphi(P_\k, f(\beta), Q_\k|\alpha) \} \in 
\Cal U^{\k}_\theta $$ 
which (again by the Lo\'s Theorem) means that
$$M\models \varphi(jP_\k, Q_\k|\alpha, j(Q_\k|\alpha))$$
In particular for $G^{*}=G\ast G_\k \ast H$ and  $G_\k\restriction \alpha \in
M[G^{*}]$ it says that
$$[G_\k\restriction \alpha]^{j} \in j(Q_\k|\alpha)_{/G^{*}} = 
j^{*}(Q_\k|\alpha).$$
\qed
\enddemo
\proclaim{Lemma 4.3}
There is a $j^*Q_\k$-generic$/M[G\ast G_\k \ast H]$ filter $H^*$
such that for every $\alpha < \kappa^+$ the condition
 $[G_\k\restriction \alpha]^j$  is in $H^*.$
\endproclaim
\demo{Proof}
Put
$$ \tilde Q = \{ q \in j^*Q_\k ;\; \forall \beta < j(\k^+): q_\beta=\emptyset
\text{ or } \max(q_\beta)\geq \k \} $$
and
$$\Cal D = \{ a \in M[G\ast G_\k \ast H];\; a\subseteq \tilde Q \text{ is 
a maximal antichain} \} .$$
It follows from the $\k^+$-c.c. of $Q_\k$ that
$$ V[G] \models \forall a \subseteq Q_\k: \text{ if }a\text{ is an antichain 
then } \exists \alpha<\k^+ :\; a\subseteq Q_\k|\alpha  $$
so
$$ M[G\ast G_\k \ast H] \models \forall a \subseteq j^* Q_\k: \text{ if }a
\text{ is an antichain 
then } \exists \alpha<j(\k^+) :\; a\subseteq j^*Q_\k|\alpha \; .$$
Moreover the cardinality of the power set of $j^*Q_\k|\alpha$ in
$M[G\ast G_\k \ast H]$ is at most $j(\k^+).$ Thus the cardinality
of $\Cal D$ in $M[G\ast G_\k \ast H]$ is $j(\k^+)$ and in $V[G\ast G_\k]$
the cardinality is $\k^+.$ Let $\langle a_\alpha;\; \alpha < \k^+ \rangle$
be an enumeration of $\Cal D$ in which each maximal antichain occurs
cofinally many times. Observe that $\tilde Q$ is $\k$-closed in $V[G\ast G_\k].$
Now it is easy to construct in $V[G\ast G_\k]$ a descending sequence of 
conditions $\langle q_\alpha ;\;\alpha<\k^+ \rangle \subseteq \tilde Q$
with the following properties:
$$(i)\; q_\alpha \in j^*(Q_\k|\alpha) ,$$
$$(ii) \; q_\alpha \leq [G_\k\restriction \alpha]^j  ,$$
$$(iii)\text{ if }a_\alpha \subseteq j^*(Q_\k|\alpha)\text{ then }
q_\alpha \text{ strenghtens a condition in }a_\alpha .$$
The sequence $\langle q_\alpha;\; \alpha<\k^+\rangle$ generates
a $j^*Q_\k$-generic$/M[G\ast G_\k \ast H]$ filter $H^*$
such that each $[G_\k\restriction \alpha]^j$ is in $H^*.$
\qed
\enddemo

It means that $p\in G_\k$ implies $j^*(p) \in H^*$ and consequently
$j^*$ is lifted to
$$ j^{**}: V[G\ast G_\k] \rightarrow M[G \ast G_\k \ast H \ast H^* ] $$
in $V[G\ast G_\k].$
We have proved that $\k$ is measurable in $V[G\ast G_\k].$

\subhead
5. Generalizations and questions
\endsubhead

We say that $S=\langle S_\lambda ; \lambda \leq \kappa \rangle$
is a generalized coherent sequence of measures if for any 
$\lambda \leq \kappa$ the set $S_\lambda$ is a set of measures over 
$\lambda$ and for any $U \in S_\lambda$
$$ j_U (S)(\lambda) = S_\lambda \restriction U = \{ V\in S_\lambda ; 
V \vartriangleleft U \} . $$
For example if each $S_\lambda$ is the set of all measures over $\lambda$ 
then the sequence is coherent. Suppose now that GCH holds, S is a generalized
coherent sequence and moreover there are separating
sets of regular cardinals $\langle X_U; U \in S_\kappa \rangle$, i.e. 
$X_U \in V$ iff $U=V$ for
$U,V \in S_\kappa .$ By a straightforward modification of our construction we
get a generic extension $V(P_{\kappa +1})$ preserving cardinalities,
cofinalities and GCH with the following properties:

1. $\langle X_U; U\in S_\kappa \rangle$ forms a maximal antichain of stationary 
subsets of $Reg$ in $V(P_{\kappa +1}) .$

2. If $U,V \in S_\kappa$ and $S\subseteq X_U , T \subseteq X_V$ are stationary
then $S<T$ (in $V(P_{\kappa+1})$) iff $U \vartriangleleft V .$ 
If $U \in S_\kappa$ and 
$S \subseteq Sing, T \subseteq X_U$ are stationary then $S<T .$
Consequently (Mitchell order)$^V(\kappa) = o^{V(P_{\kappa +1})}(\kappa) .$

Moreover if $\kappa$ has a repeat point (in a generalized sense) then
$\kappa$ is measurable in the generic extension.
It is shown in [Ba85] that any prewellordering $P$ with $|P|<\kappa$ can be 
represented as the set of all measures over $\kappa .$ That does not give us 
anything new - in that case Full Reflection again holds in the resulting model.
However recent papers of Cummings ([Cu92a], [Cu92b]) provide models with a 
rather complex structure of the Mitchell order. Using the model of [Cu92a] that 
satisfies GCH we can for example construct a generalized coherent sequence $S$ 
such that $S_\kappa$ is isomorphic to the four element poset of the type
$$\align
\circ \;\;\;\;  &\circ \\
\downarrow  \searrow &\downarrow \\
\circ \;\;\;\;  &\circ 
\endalign
 $$
Thus in the resulting model $\kappa$ is 2-Mahlo and we get two disjoint sets of   
inaccessible non-Mahlo cardinals $X_1, X_2 \subset E_0$ 
and two disjoint sets of 1-Mahlo 
cardinals $Y_1, Y_2\subset E_1$ so that for any stationary $S_1 \subseteq X_1 ,$
$S_2 \subseteq X_2$ the following holds:
$$S_1<Y_1 \text{ but } S_1 \nless Y_2 , Tr(S_1) = Y_1 \text{ (mod NS)}$$
$$S_2<Y_1 \text{ and } S_2 < Y_2 , Tr(S_2) = E_1 \text{ (mod NS)}.$$

The following question immediately comes to mind:

\proclaim{Question 1}
Does the consistency of Full Reflection at a measurable cardinal imply
the consistency of a cardinal with a repeat point?
\endproclaim

Let $V(P_{\k+1})$ be our generic extension, let $U$ denote the measure on
$\kappa$ and $C[E^{\lambda}_{\delta}]$ the filters of subsets of $\lambda$
generated in $V(P_{\k+1})$ by closed unbounded sets and the canonical 
stationary  set $E^{\lambda}_{\delta}.$ Let $\Cal F$ code all these filters, 
then
in $L[\Cal F, U]$ we get back the original measures and $U$ becomes a repeat
point of $\kappa.$ Hence we can ask a more specific question and conjecture that
the answer is yes.

\proclaim{Question 2}
Suppose that Full Reflection holds at all $\lambda \leq \kappa$, $\kappa$ is 
measurable 
and canonical stationary sets $E^\lambda _\delta$ of all orders exist. 
Let $C[E^\lambda _\delta], \Cal F , U$ be as above. Is it then true
that all $C[E^\lambda _\delta]\cap L[\Cal F, U]$ are measures in $L[\Cal F, U]$
and $U \cap L[\Cal F, U]$ is a repeat point of $\kappa$?
\endproclaim

Another way to state an equiconsistency result would be to improve our
construction so that the filters $C[E^\lambda _\delta]$ are 
$\lambda^+$-saturated. If we add this property of 
the filters to the assumptions of question 2 then using a 
method of [J84] or [JWo85] we can prove that the answer is yes. 
Unfortunately if we analyze our construction we find out that
already the filters $F^\lambda _\delta$ are not $\lambda^+$-saturated.
We can try to use the ideas of [JWo85] and instead of extensions of
$j=j_\delta : V \rightarrow M$ into $j^* : V(P_\kappa \ast Q) 
\rightarrow V(j(P_\kappa \ast Q))$ 
constructed in $V(j(P_\kappa \ast Q))$ work only with extensions 
constructed in $V(P_\kappa \ast (jP_\kappa)^\kappa) .$ We can get
the construction to work but the filters still will not be saturated. 
Hence we conjecture that the answer of the following question is no.

\proclaim{Question 3}
Is it consistent that Full Reflection holds at $\kappa$ measurable,
all canonical stationary sets $E^\kappa _\delta$ exist and the filters
$C[E^\kappa _\delta]$ are $\kappa^+$-saturated?
\endproclaim


\Refs
\ref \by {\bf [Ba85]} S.Baldwin
\paper The $\vartriangleleft$-ordering on normal ultrafilters
\jour JSL \vol 51 \yr 1985 \pages 936 -- 952 \endref

\ref \by {\bf [BTW76]} J.Baumgartner, A.Taylor and S.Wagon
\paper On splitting stationary subsets of large cardinals
\jour JSL \vol 42 \yr 1976 \pages 203--214 \endref

\ref \by {\bf [Cu92a]} J.Cummings
\paper Possible behaviors for the Mitchell ordering
\jour preliminary version \endref

\ref \by {\bf [Cu92b]} J.Cummings
\paper Possible behaviors for the Mitchell ordering II
\jour preliminary version \endref

\ref \by {\bf [HS85]} L.Harrington and S.Shelah
\paper Some exact equiconsistency results in set theory
\jour Notre Dame J.Formal logic \vol 26 \yr 1985 \pages 178--188 \endref

\ref \by {\bf [J84]} T.Jech 
\paper Stationary subsets of inaccessible cardinals
\jour Contemporary Mathematics \vol 31 \yr 1984 \pages 115--141 \endref

\ref \by {\bf [J86]} T.Jech
\book Multiple Forcing
\publ Cambridge University Press \yr 1986 \endref

\ref \by {\bf [JMa..80]} T.Jech, M.Magidor, W.Mitchell and K.Prikry
\paper Precipitous ideals
\jour JSL \vol 45 \yr 1980 \pages 1--8 \endref

\ref \by {\bf [JS90]} T.Jech, S.Shelah
\paper Full reflection of stationary sets below $\aleph_\omega$
\jour JSL \vol 55 \yr 1990 \pages 822--829 \endref

\ref \by {\bf [JS92]} T.Jech, S.Shelah
\paper Full reflection of stationary sets at regular cardinals
\jour American Journal of Mathematics, to appear \endref

\ref \by {\bf [JWo85]} T.Jech, W.H.Woodin
\paper Saturation of the closed unbounded filter on the set of 
regular cardinals
\jour Transactions of AMS \vol 292 \yr 1985 \pages 345--356 \endref

\ref \by {\bf [M82]} M.Magidor
\paper Reflecting stationary sets
\jour JSL \vol 47 \yr 1982 \pages 755--771 \endref

\ref \by {\bf [Mi74]} W.J.Mitchell
\paper Sets constructible from sequences of ultrafilters
\jour JSL \vol 39 \yr 1974 \pages 57--66 \endref

\ref \by {\bf [Mi80]} W.J.Mitchell
\book How weak is a closed unbounded filter?
\publ  Logic Colloquium '80, ed. by D.van Dalen, D.Lascar and J.Smiley, 
N.Holland 
\yr 1982 \endref

\ref \by {\bf [Mi83]} W.J.Mitchell
\paper Sets constructible from sequences of measures: revisited
\jour JSL \vol 48 \yr 1983 \pages 600--609 \endref

\ref \by {\bf [Ra82]} L.B.Radin
\paper Adding closed cofinal sequences to large cardinals
\jour Annals of Math. Logic \vol 22 \yr 1982 \pages 243--261 \endref

\ref \by {\bf [Wo92]} H.Woodin
\paper personal communication \endref

\ref \by {\bf [WoC92]} H.Woodin, J.Cummings
\book Generalised Prikry Forcings
\publ in preparation \endref


\endRefs
\enddocument
\end